\documentclass[a4paper, 12pt]{article}

\usepackage{amsfonts, amsmath, amsthm}
\usepackage{mathrsfs, latexsym}
\usepackage{bbm}

\theoremstyle{plain}
\newtheorem{thm}{Theorem}[section]

\newtheorem{lem}{Lemma}[section]

\theoremstyle{definition}

\theoremstyle{remark}
\newtheorem{rmk}{Remark}[section]

\theoremstyle{remark}
\newtheorem{exm}{Example}[section]

\numberwithin{equation}{section}
\allowdisplaybreaks

\newcommand{\N}{\mathbb{N}}

\newcommand{\R}{\mathbb{R}}
\newcommand{\C}{\mathbb{C}}
\newcommand{\Hyp}{\mathbb{H}_+}

\newcommand{\conj}[1]{\overline{#1}}
\newcommand{\Lin}{L}
\newcommand{\Non}{G}
\newcommand{\Fjc}{F_j^{\mathrm{c}}}
\newcommand{\Fjcred}{F_j^{\mathrm{c, red}}}
\newcommand{\Fcred}{F^{\mathrm{c, red}}}
\newcommand{\FjacIJK}{F_{j,a}^{\mathrm{c},I,J,K}}

\newcommand{\pa}{\partial}
\newcommand{\eps}{\varepsilon}

\newcommand{\dal}{\Box}

\DeclareMathOperator{\realpart}{\rm Re}
\DeclareMathOperator{\imagpart}{\rm Im}
\DeclareMathOperator{\supp}{\rm supp}

\begin{document}
\title{
Global existence of small amplitude solutions
to one-dimensional nonlinear Klein-Gordon systems\\
with different masses
 }

\author{
          Donghyun Kim\thanks{
              Seoul National Science Museum. 
              215, Changgyeonggung-ro, Jongno-gu, Seoul 110-360, Korea
              (E-mail: {\tt u553252d@alumni.osaka-u.ac.jp})
             }
}

\date{\today }   
\maketitle

\noindent{\bf Abstract:}\ We study the Cauchy problem for 
systems of cubic nonlinear Klein-Gordon equations 
with different masses in one space dimension.
Under a suitable structural condition on the nonlinearity, 
we will show that the solution exists globally
and decays of the rate $O(t^{-(1/2-1/p)})$
in $L^p$, $p\in[2,\infty]$ as $t$ tends to infinity
even in the case of mass resonance,
if the Cauchy data are sufficiently small, smooth
and compactly supported.
\\




\section{Introduction}  \label{sec_intro}
This paper is intended as an extension of the paper \cite{Kim}.
We consider the Cauchy problem for
a system of nonlinear Klein-Gordon equations in one space dimension
\begin{align}
 (\Box + m_j^2)u_j = F_j(u,\pa u), 
 \qquad (t,x) \in (0,\infty) \times \R,\ j=1,\ldots, N,
 \label{eq_NLKG}
\end{align}
with initial data
\begin{align}
 u_j(0,x) = \eps f_j(x), \;\; \pa_t u_j (0,x) = \eps g_j(x), 
 \qquad x \in \R,\ j=1,\ldots, N
 \label{eq_data}
\end{align}
where $\Box = \pa_t^2 - \pa_x^2$, $\pa=(\pa_t, \pa_x)$ and
$u=(u_j)_{1\le j\le N}$ is an $\R^N$-valued unknown function of 
$(t,x) \in [0,\infty)\times \R$.
Here the masses $m_1,\ldots,m_N$ are positive constants.
Without loss of generality, we assume $m_1 \le \cdots \le m_N$
throughout this paper.
The nonlinear term $F=(F_j)_{1\le j\le N}:\R^{3N}\to \R^N$ 
is assumed to be smooth and cubic around the origin, i.e., 
$$
 F_j(u,\pa u) = O((|u|+|\pa u|)^3)
 \quad \mbox{as}\ \ (u,\pa u)\to (0,0).
$$ 
For simplicity, we always suppose $f_j$, $g_j \in C_0^{\infty}(\R)$ and 
$\eps >0$ is small enough. 
This paper will focus on global existence of the solution
for \eqref{eq_NLKG}--\eqref{eq_data}
and its large-time behavior.

For $n$-space dimensional cases with 
nonlinear terms of $d$th degree
(i.e. $F_j(u, \pa u) = O((|u| + |\pa u|)^d)$ near $(u, \pa u) = (0,0)$),
it is well-known that if $d>1+2/n$,
there exists a unique global solution $u$ to the Cauchy problem
\eqref{eq_NLKG}--\eqref{eq_data}
for sufficiently small data without any restrictions 
on the nonlinearities and the masses.
This result comes from the works \cite{Kla}, \cite{Kla-Pon}, 
\cite{Sha1982} and \cite{Sha1985} done by Klainerman, Ponce and Shatah.
And it was easily verified that 
if $d>1+2/n$,
the solution is asymptotically free and satisfies
the following $L^p$-decay estimates
\begin{align}\label{decay_free}
\sum_{|\gamma|\le 1}\|\pa^\gamma u(t,\cdot)\|_{L^p} 
\le 
C\eps (1+t)^{-(n/2-1/p)}
\end{align}
for all $t\ge 0$ and $p\in[2,\infty]$ with some positive constant $C$
which does not depend on $\eps$.
Note that the condition $d>1+2/n$ can be read as $d>3$ for one-dimensional case.
When $d=3$ and $N=1$,
Yordanov \cite{Yo} introduced some kinds of nonlinearities such that the 
small data solution
for \eqref{eq_NLKG}--\eqref{eq_data} blows up in finite time.
So in general, if the nonlinear term $F_j$ is cubic near the origin,
we need some restrictions on $F_j$ to show the global existence
of the solution even if the data are sufficiently small.
As we will see below, 
it should be noted that in the case of systems, i.e. when $N\ge 2$, 
global existence and large-time behavior of the solution are affected by 
the ratio of masses as well as the structure of nonlinearities.

Let us recall some previous results briefly.
Here we concentrate our attention on the case of systems
(for the scalar case, we refer the readers to
\cite{De, Ha2008, Ka, Mo, Mo2, Su4, Yo} and the references cited therein).
Sunagawa showed in \cite{Su1} that the small solution for
\eqref{eq_NLKG}--\eqref{eq_data} exists globally 
and it approaches a free solution under some assumptions
on the nonlinearities and the masses 
(which can be seen as the {\it mass non-resonance} condition).
For example, when $N=2$ with $(m_1 - m_2)(3m_1 - m_2)\neq 0$,
the small data global existence holds for \eqref{eq_NLKG}--\eqref{eq_data}
and the solution tends to a free solution if
$F_1 = F_1(u_2, \pa u_2)$ and
$F_2 = F_2 (u_1, \pa u_1)$.
So we are more interested in the situation
that the non-resonance condition fails,
i.e. the situation that the phenomenon called {\it mass resonance} occurs 
(see \S 2 for the detail).
Normally, in the case of mass resonance, 
the situation becomes much more complicated;
for instance, according to \cite{Su1}, the solution of
\begin{align}\label{sys1}
 \left\{
 \begin{array}{l}
  (\Box+m_1^2)u_1
  =
  (\pa_t u_1)^2 (\pa_x u_2),\\
  (\Box +m_2^2)u_2
  =
  (\pa_t u_2)^2 (\pa_x u_1)
 \end{array}
 \right.
\end{align}
blows up in finite time if $m_1 = m_2$,
while the solution exists globally and tends to a free solution 
if $m_1 \neq m_2$ or $m_2 \neq 3m_1$
under some conditions on data.
To the author's knowledge, 
it is still the open problem what happens for \eqref{sys1} when $m_2 = 3m_1$.

It should be also noted that the decay estimate \eqref{decay_free}
in one space dimension is 
not a trivial one for a global small solution of \eqref{eq_NLKG}
in the mass-resonant case.
For instance, the following system 
$$
 \left\{
 \begin{array}{l}
  (\Box+m_1^2)u_1=0, \\
  (\Box +m_2^2)u_2=u_1^3
 \end{array}
 \right.
$$
admits a global solution, but according to \cite{Su2}
(see e.g., \cite{Su5} for more examples),
the second component $u_2$ does not decay 
faster than $O(t^{-1/2}\log t)$ in $L^{\infty}$
if the masses satisfy the resonant relation
(i.e., when $m_1=m_2$ or $m_2 = 3m_1$). 
So we must put some structural
conditions on the nonlinearities 
for \eqref{eq_NLKG}--\eqref{eq_data}
in order to obtain the global solution 
which satisfies \eqref{decay_free}
in the case of mass resonance.

Fang and Xue \cite{Fang} introduced
some structural conditions on the nonlinearities
called the {\it null condition}, 
under which \eqref{eq_NLKG}--\eqref{eq_data} admits a global solution
even in the mass-resonant case when $N=2$.
A pointwise asymptotic profile of the solution also can be found in \cite{Fang},
from which \eqref{decay_free} follows immediately.
But there are still some simple types of nonlinearities 
which cannot be covered by \cite{Fang};
for instance, the two component system
\begin{align}\label{sys2}
 \left\{
 \begin{array}{l}
  (\Box+m_1^2)u_1
  =
  u_1^2 u_2,\\
  (\Box +m_2^2)u_2
  =
  u_1^3,
 \end{array}
 \right.
\end{align}
with the masses satisfying the resonant relation $m_2 = 3m_1$. 
According to \cite{Su3},
the global small solution to the system
\begin{align*}
 \left\{
 \begin{array}{l}
  (\Box+m_1^2)u_1=(u_1^2+u_2^2)u_1,\\
  (\Box+m_2^2)u_2=(u_1^2+u_2^2)u_2
 \end{array}
 \right.
\end{align*}
decays like $O(t^{-1/2})$ in $L^\infty$
even in the case of mass resonance (i.e., $m_1 = m_2$),
if the data are sufficiently small, smooth and compactly supported.
Recently a structural condition on the nonlinearity is investigated in \cite{Kim}
where it is shown that 
\eqref{eq_NLKG}--\eqref{eq_data} admits a unique global solution and
it decays even strictly faster than a free solution in $L^p$, $2\le p\le \infty$
when $m_1=\cdots=m_N$.
For example, the solution of the nonlinear dissipative system
\begin{align*}
 \left\{
 \begin{array}{l}
  (\Box+m_1^2)u_1=-((\pa_t u_1)^2 + (\pa_t u_2)^2)\pa_t u_1,\\
  (\Box+m_2^2)u_2=-((\pa_t u_1)^2 + (\pa_t u_2)^2)\pa_t u_2
 \end{array}
 \right.
\end{align*}
decays rapidly of the rate $O(t^{-1/2}(\log t)^{-1/2})$
in $L^\infty$ when $m_1 = m_2$,
if the data are sufficiently small, smooth and compactly supported.

It is natural to generalize the results of \cite{Su3} or \cite{Kim} to our problem
\eqref{eq_NLKG}--\eqref{eq_data} 
which covers much more interesting situation of the mass resonance.
The aim of this paper is to introduce a new structural condition of the cubic 
nonlinearities including \eqref{sys2}
under which \eqref{eq_NLKG}--\eqref{eq_data} 
admits a unique global solution and it decays like $O\left(t^{-(1/2-1/p)}\right)$ 
in $L^p$, $2\le p \le \infty$ even in the case of mass resonance, 
if $\eps$ is small enough.

\section{Main Results} \label{sec_main}
In order to state the results, let us introduce some notations
and assumptions on the nonlinearities.
For $j\in \{1,\ldots,N\}$ and $a=(a_1, a_2, a_3)\in\{1,\ldots,N\}^3$,
we write $a\in \mathscr M_j$ if
\begin{align}\label{mass_condi}
 m_j 
 = 
 \sigma_{1} m_{a_1} 
 + \sigma_{2} m_{a_2} 
 + \sigma_{3} m_{a_3}
\end{align}
for some
$\sigma=(\sigma_1,\sigma_2,\sigma_3)\in\{+1,-1\}^3$.
Also we denote $\mathscr S_j^a$ by
$$
 \mathscr S_j^a = \left\{
 \sigma\in\{+1,-1\}^3 \;:\; m_j =  
 \sigma_{1} m_{a_1} 
 + \sigma_{2} m_{a_2} 
 + \sigma_{3} m_{a_3}
 \right\}
$$
for $a\in \mathscr M_j$.
Employing these notations, we suppose that
the nonlinear term $F_j(u, \pa_t u, \pa_x u)$ is the form of
\begin{align}\label{eq_Fj}
F_j (u, \pa_t u, \pa_x u)
=
\Fjc (u, \pa_t u, \pa_x u)
+ O\left( (|u| + |\pa_t u| + |\pa_x u|)^4\right)
\end{align}
as $(u, \pa_t u, \pa_x u) \to (0,0,0)$,
where the cubic homogeneous part of $F_j(u, \pa_t u, \pa_x u)$ 
is assumed to be
\begin{align}
\Fjc (u, \pa_t u, \pa_x u)
=
\sum_{a\in \mathscr M_j}
\sum_{|I|,|J|,|K|\le 1}
\FjacIJK(u, \pa_t u, \pa_x u)
\end{align}
and here $\FjacIJK$ takes the form of
\begin{align}\label{eq_FjabccIJK}
\FjacIJK(u, \pa_t u, \pa_x u)
=
C_{j,a}^{I,J,K} 
(\pa^I u_{a_1}) (\pa^J u_{a_2}) (\pa^K u_{a_3})
\end{align}
with some real constants $C_{j,a}^{I,J,K}$.
We may call the relation \eqref{mass_condi} with
\eqref{eq_Fj}--\eqref{eq_FjabccIJK}
{\it mass resonance} condition for one-dimensional Klein-Gordon systems
(see also \cite{Su1} or \cite{Su5} for the definition of the mass resonance).
Next, we set the upper branch of the unit hyperbola
$$
 \Hyp
 = 
 \{\omega = (\omega_0, \omega_1) \in \R^2 \; : \; 
    \omega_0 >0, \; \omega_0^2 - \omega_1^2 =1 \}.
$$
To state our conditions,
we introduce the {\it reduced nonlinearity} 
$F^{\mathrm {c, red}} = (F_j^{\mathrm {c, red}})_{1\le j\le N} :  \Hyp\times \C^N \to \C^N$
with the masses $m=(m_1,\ldots,m_N)$ by
\begin{align}\label{def_red}
 &\Fjcred (\omega, Y) \nonumber\\
 &\qquad= 
 \frac{1}{m_j}
 \sum_{a\in\mathscr M_j}
 \sum_{\sigma \in \mathscr S_j^a}
 \sum_{|I|,|J|,|K|\le 1}
 \FjacIJK
 \left(
 \mathbbm 1,
 i  \omega_0 m,
 -i \omega_1 m
 \right)
 \sigma_1^{|I|}
 \sigma_2^{|J|}
 \sigma_3^{|K|}
 Y_{a_1}^{(\sigma_{1})}
 Y_{a_2}^{(\sigma_{2})}
 Y_{a_3}^{(\sigma_{3})}
\end{align}
for $\omega \in \Hyp$ and $Y = (Y_k)_{1\le k\le N} \in \C^N$, where
$\mathbbm 1 = (1)_{1\le k\le N}$,
$i\omega_0 m = (i \omega_0 m_k)_{1\le k\le N}$,
$i\omega_1 m = (i \omega_1 m_k)_{1\le k\le N}$,
and
$Y_k^{(+1)} = Y_k$, $Y_k^{(-1)} = \overline{Y_k}$.
Also we denote by $Y\cdot Z$ 
the standard scalar product in $\C^N$ for $Y,Z \in \C^N$
and write $|Y| = \sqrt{Y\cdot Y}$ as usual.

Now we state our main results.

\begin{thm}\label{thm_SDGE}
 Assume there exists an $N \times N$ positive Hermitian matrix $A$
 such that
 \begin{align} \label{eq_decayassumption0}
  \imagpart\left(AY \cdot F^{\mathrm {c, red}} (\omega, Y)\right)
  \le 0
 \end{align}
 for all $(\omega, Y) \in \Hyp \times \C^N$.
 Then there exists $\eps_0 >0$ such that 
 if $\eps \in (0, \eps_0]$,
 the Cauchy problem \eqref{eq_NLKG}--\eqref{eq_data} admits 
 a unique global classical solution.
 Moreover, it satisfies
 \begin{align} \label{eq_decay0}
  \sum_{|\gamma|\le 1}\|\pa^\gamma u(t, \cdot)\|_{L^p} 
  \le 
  C\eps(1+t)^{-(1/2 -1/p)}
 \end{align}
 for all $t \ge 0$ and $p\in[2,\infty]$ 
 with some positive constant C which does not depend on $\eps$.
\end{thm}

\noindent We end this section by giving some examples satisfying the condition
\eqref{eq_decayassumption0}.

\begin{exm}\label{ex1}
We consider the two component system
\begin{align}\label{eq_ex1}
 \left\{
 \begin{array}{l}
  (\Box+m_1^2)u_1
  =F_1 (u, \pa_t u, \pa_x u)
  =b_1 u_1^2 u_2,\\
  (\Box +m_2^2)u_2
  =F_2 (u, \pa_t u, \pa_x u)
  =b_2 u_1^3,
 \end{array}
 \right.
\end{align}
with $m_2 = 3m_1$, $b_1, b_2 \in \R$ and $b_1 b_2 >0$.
For this system, we have
\begin{align*}
 F_1^{\mathrm{c, red}} (\omega, Y)
 =
 \frac{b_1}{m_1} \overline{Y_1}^2 Y_2,
 \qquad
 F_2^{\mathrm{c, red}} (\omega, Y)
 =
 \frac{b_2}{m_2} Y_1^3,
\end{align*}
whence  $\imagpart\left(AY \cdot F^{\mathrm {c, red}} (\omega, Y)\right)
 \le0$
for all $(\omega, Y) \in \Hyp \times \C^N$ with 
$
A = \left(\begin{matrix}
|b_2| m_1 &0\\ 0&|b_1| m_2
\end{matrix}\right)
$.\\
\end{exm}

\begin{exm}\label{ex2}
We consider the four component system
\begin{align}\label{eq_ex2}
 \left\{
 \begin{array}{l}
  (\Box+m_1^2)u_1 = F_1 (u, \pa_t u, \pa_x u) = c_1 u_2u_3u_4,\\
  (\Box +m_2^2)u_2 = F_2 (u, \pa_t u, \pa_x u) = c_2 u_3u_4u_1,\\
  (\Box +m_3^2)u_3 = F_3 (u, \pa_t u, \pa_x u) = c_3 u_4u_1u_2,\\
  (\Box +m_4^2)u_4 = F_4 (u, \pa_t u, \pa_x u) = c_4 u_1u_2u_3,
 \end{array}
 \right.
\end{align}
with $m_4 = m_1 + m_2 + m_3$.
Simply we assume that $c_1, c_2, c_3, c_4$
are positive constants.
For this system, we have 
\begin{align*}
&F_1^{\mathrm{c, red}}(\omega,Y)
=\frac{c_1}{m_1} \overline{Y_2}\overline{Y_3}Y_4,\qquad
F_2^{\mathrm{c, red}}(\omega,Y)
=\frac{c_2}{m_2} \overline{Y_3}Y_4\overline{Y_1},\\
&F_3^{\mathrm{c, red}}(\omega,Y)
=\frac{c_3}{m_3} Y_4\overline{Y_1}\overline{Y_2},\qquad
F_4^{\mathrm{c, red}}(\omega,Y)
=\frac{c_4}{m_4} Y_1 Y_2 Y_3,
\end{align*}
so that
$$
 \imagpart\left(AY \cdot \Fcred (\omega, Y)\right) \le 0
$$
for all $(\omega, Y) \in \Hyp \times \C^N$ with
$$
A = \left(\begin{matrix}
m_1/(3c_1)&0&0&0\\
0& m_2/(3c_2)&0&0\\
0&0& m_3/(3c_3)&0\\
0&0&0&m_4/c_4 
\end{matrix}\right).
$$
\end{exm}

\noindent
\newline
Hence the nonlinear terms appearing in both 
Example~\ref{ex1} and Example~\ref{ex2}
satisfy the condition \eqref{eq_decayassumption0}
with given mass relations.
Thus we can conclude that the small data global existence holds
for \eqref{eq_ex1} and \eqref{eq_ex2}.
Moreover, their solutions decay like $O(t^{-(1/2-1/p)})$ in $L^p$,
$2\le p\le \infty$ as $t\to\infty$.\\

\begin{rmk}
Note that the condition in Theorem~\ref{thm_SDGE} still holds
for \eqref{eq_ex2}
with more general coefficients $(c_k)_{1\le k\le 4}$
when $m_4 = m_1 + m_2 + m_3$.
For example, the same conclusion above is valid for
(i) $c_1, c_2, c_3 \ge 0$, $c_4>0$ and 
at least one of $c_1, c_2, c_3$ is nonzero,
(ii) $c_1 = -c_2$ and $c_3 c_4 >0$,
(iii) $c_1, c_2, c_3, c_4 <0$,
etc. with appropriately chosen $A$.
However, in the case of $c_1 = c_2 = c_3 =0$ and $c_4 \neq0$,
Theorem~\ref{thm_SDGE} fails and we know that the last component
$u_4(t)$ does not decay faster than $O(t^{-1/2}\log t)$ in $L^\infty$
according to \cite{Su2}.\\
\end{rmk}

The rest of this paper is organized as follows. 
In Section~\ref{sec_reduction}, we reduce the original problem 
\eqref{eq_NLKG}--\eqref{eq_data} by using hyperbolic coordinates.  
Section~\ref{sec_apriori} is devoted to getting a suitable a priori estimate, 
from which the small data global existence in Theorem~\ref{thm_SDGE} 
follows immediately. 
After that, we prove the time decay estimate
\eqref{eq_decay0} in Section~\ref{sec_decay}
and give some remarks in Section~\ref{sec_rmk}.
In what follows, all non-negative constants will be denoted by $C$
which may vary from line to line 
unless otherwise specified.

\section{Reduction of the Problem} \label{sec_reduction}
In this section, we perform some reduction of the problem along the idea of \cite{De} (originated from Klainerman \cite{Kla})
with a slight modification. 
In the following, we shall neglect the higher order terms of $F_j$ 
(i.e. we assume $F_j = \Fjc$) because the higher order terms do not have an essential influence on our problem
if we are interested in small amplitude solutions.

First let $B$ be a positive constant which satisfies 
$$
 \supp f_j \cup \supp g_j \subset \{ x \in \R \,:\, |x| \le B \} 
$$
and let $\tau_0 > 1 + 2B$. We start with the fact that we may treat the 
problem as if the Cauchy data are given on the upper branch of the hyperbola
$$
 \{ (t,x) \in \R^2 \,:\, (t+2B)^2 - x^2 = \tau_0^2, \; t >0 \}
$$
and they are sufficiently smooth, small, compactly-supported. 
This is a consequence of the classical local existence theorem 
and the finite speed of propagation 
(see e.g., Proposition 1.4 of \cite{De} for the detail). 
Next, let us introduce the hyperbolic coordinates
$(\tau, z ) \in [\tau_0 , \infty) \times \R$ in the interior of the light cone, i.e.,
$$
 t+2B = \tau \omega_0 (z),\quad x = \tau \omega_1 (z)
$$ 
for $|x| < t + 2B$ where
$\omega(z) = (\omega_0(z), \omega_1(z))=(\cosh z, \sinh z) \in \Hyp$. 
Note that $\omega_k (z)$ satisfies
$$
 |\omega_k (z)| \le Ce^{|z|}
$$
for all $z \in \R$ and $k=0,1$.
Then we can easily check that
$$
 \pa_t = \omega_0(z) \pa_\tau - \frac{1}{\tau}\omega_1(z) \pa_z,\quad 
 \pa_x = -\omega_1(z) \pa_\tau + \frac{1}{\tau}\omega_0(z) \pa_z,\quad
 \dal = \pa_\tau^2 + \frac{1}{\tau} \pa_\tau - \frac{1}{\tau^2} \pa_z^2.
$$
We also take a weight function $\chi(z) \in C^\infty(\R)$ satisfying
$$
 0 < \chi (z) \le C_0 e^{-\kappa |z|} \;\;\text{and}\;\; 
 | \pa_z^k \chi (z) | \le C_k \chi (z) \;\;(k=1,2,\cdots)
$$
with a large parameter $\kappa \gg 1$ and positive constants $C_0, C_j$.
With this weight function, 
let us define the new unknown function $v(\tau, z) = (v_j(\tau,z))_{1\le j\le N}$ by
\begin{equation*}
 u_j (t,x)
 =
 \frac{\chi(z)}{\sqrt \tau} v_j (\tau,z).
\end{equation*}
Then we see that $v_j$ satisfies 
$\Lin_j v_j=\Non_j$ 
if $u$ solves \eqref{eq_NLKG}, where
\begin{align}
 \Lin_j
 = 
 \pa_\tau^2 - \frac{1}{\tau^2} 
 \left( 
   \pa_z^2 + 2 \frac{\chi' (z)}{\chi(z)} \pa_z + \frac{\chi''(z)}{\chi(z)} 
   - \frac{1}{4} 
 \right) 
 +m_j^2
 \label{eq_op_lin}
\end{align}
and
$$
 \Non_j 
 = 
 \frac{\chi(z)^2}{\tau}
 \Fjc(v,\, \omega_0(z)\pa_\tau v,\, -\omega_1(z)\pa_\tau v)
 + {Q}_j. 
$$
Here ${Q}_j$ takes the form of
$$
Q_j
 =
 \sum_{k_1, k_2, k_3=1}^{N}\ 
 \sum_{\nu_1, \nu_2=0}^{1}\ 
 \sum_{\substack{0\le \eta_1\le \nu_1\\
                 0\le \eta_2\le \nu_2\\
                 0\le \eta_3\le 1}}
 \frac{q_{j, k_1, k_2, k_3}^{\nu_1,\nu_2,\eta_1,\eta_2,\eta_3}(z)}
      {\tau^{2+\nu_1+\nu_2}}
 \bigl(\partial_{\tau}^{1-\nu_1} \partial_z^{\eta_1} v_{k_1}\bigr)
 \bigl(\partial_{\tau}^{1-\nu_2} \partial_z^{\eta_2} v_{k_2}\bigr)
 \bigl(\partial_z^{\eta_3} v_{k_3}\bigr)
$$
with some 
$q_{j,k_1,k_2,k_3}^{\nu_1,\nu_2,\eta_1,\eta_2,\eta_3} \in C^{\infty}(\R)$ 
satisfying 
$$
 |\pa_z^l q_{j,k_1,k_2,k_3}^{\nu_1,\nu_2,\eta_1,\eta_2,\eta_3} (z)|\le 
 C_l e^{(3-2\kappa)|z|}
$$
for $l\in\{0,1,2,\ldots\}$.
As we shall see in \eqref{est_Qj} below, 
$Q_j$ can be regarded as a remainder, 
while the first term of $G_j$ plays a role as a main term.

At last, the original problem \eqref{eq_NLKG}--\eqref{eq_data} is reduced to
\begin{align}
 \left\{
 \begin{array}{l}
  \Lin_j v_j 
  = 
  \Non_j (\tau, z, v, \pa_\tau v, \pa_z v), \quad 
  \tau > \tau_0, \; z \in \R, \; j=1, \ldots, N, 
  \\
  (v_j, \pa_\tau v_j) |_{\tau=\tau_0}
  =
  (\eps \tilde{f_j}, \eps \tilde{g_j}), \quad 
  z \in \R, \; j=1, \ldots, N 
 \end{array}
 \right.
\label{eq_reduced}
\end{align}
where $\tilde{f_j}, \tilde{g_j}$ are sufficiently smooth functions of $z$ 
with compact support.

\section{A Priori Estimate} \label{sec_apriori}
This section is devoted to getting an a priori estimate 
for the solution of the reduced problem \eqref{eq_reduced}
under the condition \eqref{eq_decayassumption0}.
As in \cite{Su4}, \cite{Kawa} and \cite{Kim}, we set
$$
 M(T)
 =
 \sup_{(\tau,z)\in[\tau_0,T)\times\R} 
 \left( |v(\tau,z)| + |\pa_\tau v(\tau,z)| 
 + \frac{1}{\tau} |\pa_z v(\tau,z)| \right)
$$
for the smooth solution $v(\tau,z)$ to \eqref{eq_reduced} 
on $\tau \in [\tau_0, T)$. We will prove the following:

\begin{lem} \label{lem_apriori}
 Under the assumption of Theorem~\ref{thm_SDGE}, 
 there exist $\eps_1 >0$ and $C^* >0$ such that 
 $M(T) \le \sqrt \eps$ implies $M(T) \le C^* \eps$ for any $\eps \in (0, \eps_1]$.
 Here $C^*$ is independent of $T$.
\end{lem}

Once this lemma is proved, 
we can derive the global existence part of Theorem~\ref{thm_SDGE} 
in the following way: By taking $\eps_0 \in (0, \eps_1]$ so that 
$2C^* \eps_0^{1/2} \le 1$, we deduce that $M(T) \le \eps^{1/2}$ implies 
$M(T) \le \eps^{1/2} /2$ for any $\eps \in (0,\eps_0]$. 
Then by the continuity argument, 
we have $M(T) \le C^* \eps$ as long as the solution exists. 
Therefore the local solution to \eqref{eq_reduced} can be extended to the global one.
Going back to the original variables, 
we deduce the small data global existence for \eqref{eq_NLKG}--\eqref{eq_data}.\\

To prove Lemma~\ref{lem_apriori}, we introduce some lemmas here.

\begin{lem} \label{lem_energy}
Assume $M(T)\le \sqrt \eps$ and let $s\ge 3$, $0<\delta<1/3$.
Then the estimate
\begin{align*}
 E_s (\tau)
 \le
 C\eps^2 \tau^\delta
\end{align*}
holds for $\tau \in [\tau_0, T)$, where the energy is defined by
$$
 E_s (\tau) 
 = 
 \sum_{k=0}^s \sum_{j=1}^N \frac{1}{2} \int_{\R} 
 | \pa_\tau \pa_z^k v_j(\tau, z) |^2 
 + \left| \frac{\pa_z}{\tau} \pa_z^k v_j(\tau, z) \right|^2 
 + m_j^2 |\pa_z^k v_j(\tau, z) |^2 \, dz.
$$
\end{lem}

\noindent We omit the proof of Lemma~\ref{lem_energy} here
because it is almost same as that of the previous works 
(\cite{Su3}, \cite{Su4} or \cite{Kim}, see also Appendix).\qed\\

\noindent For the next step, we introduce the $\C^N$-valued function 
$\alpha = (\alpha_j)_{1\le j\le N}$ by
\begin{align}\label{def_alpha}
 \alpha_j (\tau, z) 
 = 
 e^{-i m_j \tau} \left(1 -\frac{i \pa_\tau}{m_j} \right) v_j (\tau, z), 
 \qquad (\tau, z) \in [\tau_0, T) \times \R, \quad j=1, \ldots, N
\end{align}
for the solution $v(\tau,z)$ to \eqref{eq_reduced}.

\begin{lem} \label{lem_alpha}
Under the assumption of Lemma~\ref{lem_energy}, we have
$$
 \left| \frac{\pa\alpha}{\pa\tau}(\tau,z)\right|
 \le
  \frac{C\eps}{\tau}
$$
for $(\tau,z)\in[\tau_0,T)\times\R$, where $\alpha$ is given by \eqref{def_alpha}.
\end{lem}

\noindent{\it Proof.}
First we note that
\begin{align} \label{eq_alphaDE1}
 \frac{\pa\alpha_j}{\pa\tau}
 &= 
 -\frac{ie^{-i m_j \tau}}{m_j} (\pa_\tau^2 +m_j^2) v_j \nonumber \\
 &=
 -\frac{ie^{-i m_j \tau}}{m_j} 
 \left( 
  \Non_j (\tau,z,v,\pa_\tau v,\pa_z v) 
  + \frac{1}{\tau^2} 
  \left( 
    \pa_z^2 + 2 \frac{\chi'(z)}{\chi(z)} \pa_z 
    + \frac{\chi''(z)}{\chi(z)} - \frac{1}{4} 
  \right) v_j 
 \right) \nonumber \\
 &=
 -\frac{ie^{-i m_j \tau}}{m_j} \frac{\chi(z)^2}{\tau} 
 \Fjc (v, \omega_0 (z) \pa_\tau v, - \omega_1 (z) \pa_\tau v) 
 + \frac{R_j (\tau, z)}{\tau^2},
\end{align}
where
\begin{align}\label{def_Rj}
 R_j (\tau, z) 
 &=
 -\frac{ie^{-im_j\tau}}{m_j} \tau^2 Q_j
 -\frac{ie^{-im_j\tau}}{m_j} \left( 
     \pa_z^2 + 2 \frac{\chi'(z)}{\chi(z)} \pa_z +\frac{\chi''(z)}{\chi(z)} 
     - \frac{1}{4} 
   \right) v_j.
\end{align}
Using the assumption $M(T) \le \sqrt \eps$, we have 
$$
 \frac{1}{m_j}\left| 
 \frac{\chi(z)^2}{\tau}\Fjc(v, \omega_0(z) \pa_\tau v, -\omega_1(z) \pa_\tau v) 
 \right| 
 \le
 \frac{C}{\tau} M(T)^3 
 \le 
 \frac{C}{\tau}\eps^{3/2}.
$$
The standard Sobolev embedding and Lemma~\ref{lem_energy} yields
\begin{align}\label{est_Qj}
 \frac{1}{m_j}|Q_j|
 \le
 \frac{C}{\tau^2}(|\pa_\tau v|^3 + |\pa_z v|^3 +|v|^3)
 \le
 \frac{C}{\tau^2} E_2(\tau)^{3/2}
 \le 
 C\eps^3 \tau^{3\delta/2-2}
\end{align}
and 
$$
 \frac{1}{m_j}\left| \left( 
  \pa_z^2 + 
  \frac{2 \chi'(z)}{\chi(z)}\pa_z
  +\frac{\chi''(z)}{\chi(z)}-\frac{1}{4} 
 \right) v_j \right| 
 \le
 C \sum_{k=0}^2 | \pa_z^k v | 
 \le
 C E_3(\tau)^{1/2} 
 \le 
 C \eps \tau^{\delta/2}
$$
so that
\begin{align} \label{eq_remainder}
 |R_j (\tau, z) | 
 \le
 C \eps \tau^{3\delta/2}.
\end{align}
Combining all together, we obtain
\begin{align*}
 \left| \frac{\pa \alpha_j}{\pa \tau} (\tau, z) \right| 
 \le
 \frac{C\eps^{3/2}}{\tau} + \frac{C\eps}{\tau^{2-3\delta/2}}
 \le 
 \frac{C\eps}{\tau},
\end{align*}
which proves Lemma~\ref{lem_alpha}.\qed\\

Next, we introduce the lemma which will be used
to estimate some non-resonant terms appearing  in \eqref{eq_alphaAalpha} below.
\begin{lem} \label{lem_nonresonant}
Let the function $\alpha$ be given in \eqref{def_alpha} and assume
$M(T)\le \sqrt\eps$. Then for any $k_1, \ldots, k_4 \in \{1, \ldots, N\}$,
$\beta_1, \ldots, \beta_4 \in \{+1,-1\}$ and
$b \in \R\backslash\{0\}$, we have the estimate
$$
 \left|
 \int_{\tau_0}^\tau
  \left(\prod_{l=1}^4 \alpha_{k_l}^{(\beta_l)}(\tilde\tau,z)\right)
  \frac{e^{i{b}\tilde\tau}}{\tilde\tau}\,d\tilde\tau\right|
 \le
 C\eps^2
$$
where 
$\alpha_{k_l}^{(+1)}=\alpha_{k_l}$ and
$\alpha_{k_l}^{(-1)}=\conj{\alpha_{k_l}}$.
\end{lem}
\noindent{\it Proof.}
We observe that
\begin{align*}
 \int_{\tau_0}^\tau
  \left(\prod_{l=1}^4 \alpha_{k_l}^{(\beta_l)}(\tilde\tau,z)\right)
  \frac{e^{ib\tilde\tau}}{\tilde\tau}\,d\tilde\tau 
 &=
 \int_{\tau_0}^\tau (\pa_\tau U_1)(\tilde\tau,z) + U_2(\tilde\tau,z)\,d \tilde\tau \\
 &=
 U_1 (\tau,z) - U_1 (\tau_0,z) + \int_{\tau_0}^\tau U_2 (\tilde\tau,z)\,d\tilde\tau,
\end{align*}
where 
\begin{align*}
 &U_1(\tau,z) 
 = 
 \left(\prod_{l=1}^4 \alpha_{k_l}^{(\beta_l)}(\tau,z)\right)
  \frac{e^{i b \tau}}{ i b \tau },\\ 
 &U_2(\tau,z) 
 = 
 \left(\prod_{l=1}^4 \alpha_{k_l}^{(\beta_l)}(\tau,z)\right)
  \frac{e^{i b \tau}}{i b \tau^2} 
 -
  \sum_{l=1}^4 \Bigg{(} \left(\pa_\tau \alpha_{k_l}^{(\beta_l)}(\tau,z)\right)
   \prod_{n \in\{1,2,3,4\} \backslash \{l\}}
   \alpha_{k_n}^{(\beta_n)}(\tau,z)\Bigg{)}
   \frac{e^{i b \tau}}{i b\tau}.
\end{align*}
Using $M(T)\le\sqrt\eps$ and Lemma~\ref{lem_alpha}, we have 
$$
 |U_1(\tau,z)| 
 \le 
 \frac{C\eps^2}{\tau}, 
 \qquad |U_2(\tau,z)|
 \le 
 \frac{C\eps^2}{\tau^2} 
  +\frac{C\eps}{\tau} \cdot \eps^{3/2} \cdot \frac{1}{\tau} 
 \le 
 \frac{C\eps^2}{\tau^2}
$$
and that
$$
 \left|\int_{\tau_0}^\tau
  \left(\prod_{l=1}^4 \alpha_{k_l}^{(\beta_l)}(\tilde\tau,z)\right)
  \frac{e^{i b\tilde\tau}}{\tilde\tau}\,d\tilde\tau\right|
 \le
 C\eps^2 \left( 1+ \int_1^\infty \frac{d \tilde\tau}{\tilde\tau^2} \right),\\
$$
which proves Lemma~\ref{lem_nonresonant}.\qed\\

In the next lemma, we will derive a certain system of ordinary differential
equations for $\alpha(\tau,z)$ defined by \eqref{def_alpha} 
regarding $z$ as a parameter.
\begin{lem} \label{lem_dalphadtau}
For each $j\in\{1,\ldots,N\}$, we have
$$
 \frac{\pa\alpha_j}{\pa\tau}(\tau,z)
 =
  -\frac{i\chi(z)^2}{8\tau}\Fjcred(\omega(z),\alpha(\tau,z))
  +S_j(\tau,z)+\frac{R_j(\tau,z)}{\tau^2}
$$
for $(\tau,z)\in [\tau_0,T)\times\R$ with 
$\Fjcred$ and $R_j$ given by 
\eqref{def_red} and \eqref{def_Rj} respectively. Here the
non-resonant term $S_j$ is the form of
\begin{align*}
 &S_j (\tau, z)\\
 &=
 -\frac{i\chi(z)^2}{m_j\tau}
 \sum_{a\in\mathscr M_j}
 \sum_{(\varsigma_{1}, \varsigma_{2}, \varsigma_{3}) \in \mathscr T_j^a}
 \sum_{|I|,|J|,|K|\le 1}
 \Omega_{j,a}^{I,J,K}(z) 
 \varsigma_1^{|I|} 
 \varsigma_2^{|J|} 
 \varsigma_3^{|K|}
 \alpha_{a_1}^{(\varsigma_1)}
 \alpha_{a_2}^{(\varsigma_2)}
 \alpha_{a_3}^{(\varsigma_3)}
 e^{i( \varsigma_1 m_{a_1} + \varsigma_2 m_{a_2}
 +\varsigma_3 m_{a_3} -m_j)\tau}
\end{align*}
where $\mathscr T_j^a=\{+1,-1\}^3 \setminus \mathscr S_j^a$
and the coefficients are given by
$$
 \Omega_{j,a}^{I,J,K}(z)
 =
 \frac{1}{8}\FjacIJK\left(
 \mathbbm 1, i\omega_0(z) m, -i\omega_1(z) m\right).
$$
\end{lem}
\noindent{\it Proof.}
For the proof, 
we introduce a new function $H_j(\theta,z)$ which is defined by 
$$
 H_j(\theta,z) 
 = 
 e^{-i m_j \theta} \Fjc (
 \realpart (\alpha e^{i m\theta}),\, 
 -\omega_0(z) \imagpart(\alpha m e^{i m\theta}),\, 
 \omega_1(z) \imagpart(\alpha m e^{i m\theta}))
$$
with $\alpha$ being regarded as a parameter for the moment.
Here we denote
$\realpart (\alpha e^{i m\theta})$
=
$((\realpart(\alpha_k e^{i m_k \theta}))_{1\le k\le N}$,
$\realpart (\alpha m e^{i m\theta})$
=
$((\realpart(\alpha_k m_k e^{i m_k \theta}))_{1\le k\le N}$
and
$\imagpart (\alpha m e^{i m\theta})$
=
$((\imagpart(\alpha_k m_k$ $e^{i m_k \theta}))_{1\le k\le N}$
respectively.
Then from \eqref{eq_alphaDE1}, we can easily check that
\begin{align}\label{eq_alphaDEH}
 \frac{\pa\alpha_j}{\pa\tau}(\tau,z)
 =
 -\frac{i\chi(z)^2}{m_j\tau} H_j(\tau,z) + \frac{R_j(\tau,z)}{\tau^2}.
\end{align}
From now on, we will decompose $H_j$ into
the resonant terms and the nonresonant terms.
Since $\Fjc$ is cubic, we may expand $H_j$ as
\begin{align*}
 &H_j(\theta,z)\\
 &\quad=
 \sum_{a\in\mathscr M_j}
 \sum_{|I|,|J|,|K|\le 1}
 \FjacIJK\left(
 \realpart (\alpha e^{i m\theta}),\, 
 -\omega_0(z) \imagpart(\alpha me^{i m\theta}),\, 
 \omega_1(z) \imagpart(\alpha me^{i m\theta})\right)
 e^{-im_j\theta}\nonumber\\
 &\quad= 
 \sum_{a\in\mathscr M_j}
 \sum_{\varsigma_1, \varsigma_2, \varsigma_3 \in \{+1, -1\}}
 \sum_{|I|,|J|,|K|\le 1}
 \Omega_{j,a}^{I,J,K}(z) 
 \Psi_1 \alpha_{a_1}^{(\varsigma_1)}
 \Psi_2 \alpha_{a_2}^{(\varsigma_2)}
 \Psi_3 \alpha_{a_3}^{(\varsigma_3)}
 e^{i(\varsigma_1 m_{a_1}+ 
 \varsigma_2 m_{a_2}+\varsigma_3 m_{a_3})\theta}
 e^{-im_j\theta}
\end{align*}
where $\alpha_k^{(+1)} = \alpha_k$, 
$\alpha_k^{(-1)} = \overline{\alpha_k},$
and
$\Psi_1 = \Psi_1(\varsigma_1, |I|)$
which corrects the sign of $\alpha_{a_1}$ given by
\begin{align*}
 \Psi_1(\varsigma_1, |I|) = \left\{
 \begin{array}{l}
   +1\qquad \text{if }(\varsigma_1, |I|) = (+1, 0),\\
   +1\qquad \text{if }(\varsigma_1, |I|) = (+1, 1),\\
   +1\qquad \text{if }(\varsigma_1, |I|) = (-1, 0),\\
   -1\qquad \text{if }(\varsigma_1, |I|) = (-1, 1),
 \end{array}
 \right.
\end{align*}
or simply $\Psi_1 = \varsigma_1^{|I|}$.
Similarly $\Psi_2 = \varsigma_2^{|J|}$ and $\Psi_3 = \varsigma_3^{|K|}$.
Remembering $a\in\mathscr M_j$, we can rewrite $H_j$ as
\begin{align*}
 &H_j(\theta,z)\\
 &\quad=
 \sum_{a\in\mathscr M_j}
 \sum_{\sigma\in\mathscr S_j^a}
 \sum_{|I|,|J|,|K|\le 1}
 \Omega_{j,a}^{I,J,K}(z) 
 \sigma_1^{|I|} 
 \sigma_2^{|J|} 
 \sigma_3^{|K|}
 \alpha_{a_1}^{(\sigma_{1})}
 \alpha_{a_2}^{(\sigma_{2})}
 \alpha_{a_3}^{(\sigma_{3})}\\
 &\qquad+
 \sum_{a\in\mathscr M_j}
 \sum_{(\varsigma_{1}, \varsigma_{2}, \varsigma_{3}) \in 
 \mathscr T_j^a}
 \sum_{|I|,|J|,|K|\le 1}
 \Omega_{j,a}^{I,J,K}\left(z\right) 
 \varsigma_1^{|I|} 
 \varsigma_2^{|J|} 
 \varsigma_3^{|K|}
 \alpha_{a_1}^{(\varsigma_1)}
 \alpha_{a_2}^{(\varsigma_2)}
 \alpha_{a_3}^{(\varsigma_3)}
 e^{i(\varsigma_1 m_{a_1} + \varsigma_2 m_{a_2} + \varsigma_3 m_{a_3}
 -m_j)\theta}.
\end{align*}
Note that $\varsigma_1 m_{a_1}+ \varsigma_2 m_{a_2}
+ \varsigma_3 m_{a_3} - m_j \neq 0$
for $(\varsigma_1, \varsigma_2, \varsigma_3)\in\mathscr T_j^a$.
Observing
$$
 \sum_{a\in\mathscr M_j}
 \sum_{\sigma\in\mathscr S_j^a}
 \sum_{|I|,|J|,|K|\le 1}
 \Omega_{j,a}^{I,J,K}(z) 
 \sigma_1^{|I|} 
 \sigma_2^{|J|} 
 \sigma_3^{|K|}
 \alpha_{a_1}^{(\sigma_{1})}
 \alpha_{a_2}^{(\sigma_{2})}
 \alpha_{a_3}^{(\sigma_{3})} 
 =
 \frac{m_j}{8} \Fjcred(\omega(z),\alpha)
$$
and with \eqref{eq_alphaDEH}, we arrive at Lemma~\ref{lem_dalphadtau}.
\qed\\

\noindent
\textbf{Proof of Lemma~\ref{lem_apriori}.}
Now we are ready to prove Lemma~\ref{lem_apriori}.
First, with the matrix $A$ in Theorem~\ref{thm_SDGE}, we note that
\begin{align} \label{eq_eigen}
 |Y\cdot AZ|^2
 \le 
 (Y\cdot AY)(Z\cdot AZ), 
 \qquad 
 \lambda_* |Y|^2
 \le 
 Y\cdot AY
 \le 
 \lambda^*|Y|^2
\end{align}
holds for any $Y, Z \in \C^N$, where $\lambda^*$ (resp. $\lambda_*$) is the 
largest (resp. smallest) eigenvalue of $A$.
Using the notations 
$S = (S_j)_{1\le j\le N}$, $R = (R_j)_{1\le j\le N}$, 
it follows from 
Lemma~\ref{lem_dalphadtau},
\eqref{eq_decayassumption0},
\eqref{eq_eigen}
and
\eqref{eq_remainder}
that
\begin{align*}
 \frac{\pa}{\pa \tau} \left( \alpha(\tau,z)\cdot A\alpha(\tau,z) \right)
 &= 
 2\realpart \left(\pa_\tau \alpha(\tau,z) \cdot A\alpha(\tau,z)\right) \\
 &=
 \frac{\chi(z)^2}{4\tau}
 \imagpart\left(\Fcred(\omega(z),\alpha)\cdot A\alpha\right)
  +2\realpart \left( S\cdot A\alpha\right)
  +\frac{2}{\tau^2} \realpart \left(R\cdot A\alpha\right) \\
 &\le
  2\realpart \left( S\cdot A\alpha\right)
  +\frac{2}{\tau^2} \left|R\cdot A\alpha\right| \\
 &\le
 2\realpart \left( S\cdot A\alpha\right)
  +\frac{1}{\tau^2} \left( \alpha \cdot A\alpha + R\cdot AR  \right)\\
 &\le
 2\realpart \left( S\cdot A\alpha\right) 
  +\frac{1}{\tau^2} \left(\alpha\cdot A\alpha\right)
  +\frac{C \eps^2}{\tau^{2- 3\delta}}
\end{align*}
which yields
\begin{align}\label{eq_alphaAalpha}
 \alpha(\tau,z) \cdot A\alpha(\tau,z)
 &\le
 C\eps^2
  +
 2 \left| \int_{\tau_0}^\tau \realpart \left( S\cdot A\alpha\right)\,d\tilde\tau \right| 
  +
 \int_{\tau_0}^\tau \alpha(\tilde\tau,z)\cdot A\alpha(\tilde\tau,z)\,
 \frac{d\tilde{\tau}}{\tilde\tau^2}\\
 &\le
 C\eps^2
  +
 \int_{\tau_0}^\tau \alpha(\tilde\tau,z)\cdot A\alpha(\tilde\tau,z)\,
 \frac{d\tilde{\tau}}{\tilde\tau^2}\nonumber
\end{align}
for $\tau \in [\tau_0, T)$, where we used Lemma~\ref{lem_nonresonant}
to estimate the second term of the right-hand side in \eqref{eq_alphaAalpha}.
Thus the Gronwall lemma and \eqref{eq_eigen} yields 
\begin{align}\label{eq_est_alpha}
 \sup_{(\tau,z)\in[\tau_0,T)\times\R}|\alpha(\tau,z)|
 \le
 C\eps.
\end{align}
In view of the relations
$$
 \left(|v(\tau,z)|^2 + |\pa_\tau v(\tau,z)|^2 \right)^{1/2}
 \le
 C|\alpha(\tau,z)|
$$
and
\begin{align}
 \frac{1}{\tau}|\pa_z v(\tau,z)|
 \le
 \frac{C}{\tau}E_2(\tau)^{1/2}
 \le
 C\eps \tau^{-(1-\delta/2)}
 \le
 C\eps,
 \label{eq_est_pa_z_v}
\end{align}
we reach Lemma~\ref{lem_apriori}.\qed

\section{Proof of the Decay Estimates}  \label{sec_decay}
Now we are in a position to prove the time decay estimate
\eqref{eq_decay0} under the condition \eqref{eq_decayassumption0}.
First, we remember that our change of variable is
\begin{align} \label{eq_substitution}
 u_j(t,x)
 =
  \frac{\chi(z)}{\sqrt{\tau}}v_j(\tau,z)
 =
 \frac{\chi(z)}{\sqrt\tau} \realpart(\alpha_j(\tau,z)e^{im_j\tau})
\end{align}
with
$t+2B = \tau\omega_0(z)$,
$x=\tau\omega_1(z)$
for $|x| < t+2B$, 
and that $u(t,\cdot)$ is supported on $\{x \in \R :|x|\le t+B\}$.
So it follows from \eqref{eq_est_alpha} and \eqref{eq_substitution} that
\begin{align}\label{eq_udecay0}
 |u(t,x)|
 \le
 \frac{\chi(z)|\alpha(\tau,z)| \omega_0(z)^{1/2}}{\sqrt{\tau}\omega_0(z)^{1/2}}
 \le
 C\eps(1+t)^{-1/2}.
\end{align}
Using \eqref{eq_udecay0} and the finite propagation speed, we have 
\begin{align*}
 \|u(t, \cdot) \|_{L^p(\R)} 
 &=
 \|u(t, \cdot) \|_{L^p(\{ |x| \le t+B\})} \\
 &\le
 C\eps (1+t)^{-1/2}
  \left( \int_{|x| \le t+B}1\,dx \right)^{1/p}\\
 &\le
 C\eps (1+t)^{-(1/2-1/p)}
\end{align*}
for $p \in [2,\infty]$.
Also we note that $\pa_t u_j$ and $\pa_x u_j$ can be written as
\begin{align*}
 \pa_t u_j(t,x)
 &=
 \left( \omega_0(z)\pa_\tau - \frac{1}{\tau}\omega_1(z)\pa_z\right)
 \left(\frac{\chi(z)}{\sqrt{\tau}}v_j(\tau,z) \right)\\
 &=
 -\frac{\chi(z)\omega_0(z)^{3/2}}{\sqrt{\tau\omega_0(z)}}
   \imagpart(\alpha_j m_j e^{im_j\tau})
 -\frac{\chi(z)\omega_0(z)}{\tau^{3/2}}
  \left( \frac{\omega_1(z)}{\omega_0(z)}\pa_z
  +\frac{\omega_1(z)}{\omega_0(z)} \frac{\chi'(z)}{\chi(z)}
  +\frac{1}{2} \right) v_j
\end{align*}
and
\begin{align*}
 \pa_x u_j(t,x)
 &=
 \left( -\omega_1(z)\pa_\tau + \frac{1}{\tau}\omega_0(z)\pa_z\right)
 \left(\frac{\chi(z)}{\sqrt{\tau}}v_j(\tau,z) \right)\\
 &=
 \frac{\frac{\omega_1(z)}{\omega_0(z)}\chi(z)\omega_0(z)^{3/2}}{\sqrt{\tau\omega_0(z)}}
   \imagpart(\alpha_j m_j e^{im_j\tau})
 + \frac{\chi(z)\omega_0(z)}{\tau^{3/2}}
  \left( \pa_z
  +\frac{\chi'(z)}{\chi(z)}
  + \frac{1}{2} \frac{\omega_1(z)}{\omega_0(z)} \right) v_j.
\end{align*}
Thus using \eqref{eq_est_pa_z_v} and \eqref{eq_est_alpha}, we get
\begin{align*}
 \sum_{|\gamma|=1} \left| \pa^\gamma u(t,x) \right|
 &\le
 \frac{C\eps}{\sqrt{\tau\omega_0(z) }}
 +
 \frac{C\eps\chi(z)\omega_0(z)^{1+ 3/2-\delta}}
      {(\tau\omega_0(z))^{3/2-\delta}}
 +
 \frac{C\eps\chi(z)\omega_0(z)^{1+3/2}}
      {(\tau\omega_0(z))^{3/2}}\\
 &\le
 C\eps (1+t)^{-1/2}
 +C\eps(1+t)^{-3/2+\delta}
 +C\eps(1+t)^{-3/2}\\
 &\le
 C\eps (1+t)^{-1/2}
\end{align*}
and from the finite propagation speed, we finally obtain the desired estimate
$$
 \sum_{|\gamma|\le 1} \| \pa^\gamma u(t, \cdot) \|_{L^p(\R)}
 \le
 C \eps (1+t)^{-(1/2-1/p)}
$$
for $p \in [2,\infty]$.
This completes the proof of Theorem~\ref{thm_SDGE}.\qed\\

\section{Further Remarks}  \label{sec_rmk}
We may modify Theorem 2.2 and Theorem 2.3 in \cite{Kim} to our problem
\eqref{eq_NLKG}--\eqref{eq_data}.
Actually, by the similar method of \cite{Kim},
we can prove the following theorem:
\begin{thm}\label{thm_decay}
 Let $\Fcred$ be given by \eqref{def_red}. Assume there exist an $N \times N$
 positive Hermitian matrix $A$ and a positive constant $\widetilde C$ such that
 \begin{align} \label{eq_decay_a1}
  \imagpart\left(AY \cdot F^{\mathrm {c, red}} (\omega, Y)\right)
  \le - \widetilde C \omega_0^k |Y|^4
  \qquad
  \forall (\omega, Y) \in \Hyp \times \C^N
 \end{align}
 for some $k\in\{1,3\}$.
 Then the global solution of
 \eqref{eq_NLKG}--\eqref{eq_data} satisfies
 \begin{align*}
  \sum_{|\gamma|\le \left\lfloor \frac{k}{2} \right\rfloor}
  \|\pa^\gamma u(t, \cdot)\|_{L^p} 
  \le 
  C \frac{(1+t)^{-(1/2 -1/p)}}{\sqrt{\log (2+t)}}
 \end{align*}
 for all $t \ge 0$ and $p\in[2,\infty]$.
\end{thm}

We can find some nonlinearity $F$ which satisfy the condition
\eqref{eq_decay_a1} under the mass resonance
(except the case of $m_1 = \cdots = m_N$ which can be covered by \cite{Kim});
for instance, we may set $F_j$ to contain some nonlinear dissipative terms like
\begin{align}\label{eq_dissipation}
 \left\{
 \begin{array}{l}
  (\Box+m_1^2)u_1
  =F_1(u,\pa_t u,\pa_x u)
  =-|\pa_t u|^2 \pa_t u_1
    -(\pa_t u_1)^2 \pa_t u_2\\
  (\Box+m_2^2)u_2=
  F_2(u,\pa_t u,\pa_x u)=
  -|\pa_t u|^2 \pa_t u_2
    +(\pa_t u_1)^3\\
 \end{array}
 \right.
\end{align}
with $m_2 = 3m_1$ whose reduced nonlinearities are
\begin{align*}
 &F_1^{\mathrm{c, red}}(\omega,Y)
 =
 -3i\omega_0^3 m_1^2 |Y_1|^2 Y_1
 -2i\omega_0^3 m_2^2 |Y_2|^2 Y_1
 +i\omega_0^3 m_1 m_2 \overline{Y_1}^2 Y_2,\\
 &F_2^{\mathrm{c, red}}(\omega,Y)
 =
 -2i\omega_0^3 m_1^2 |Y_1|^2 Y_2
 -3i\omega_0^3 m_2^2 |Y_2|^2 Y_2
 -i \omega_0^3 \frac{m_1^3}{m_2} Y_1^3.
\end{align*}
Then we see that
\begin{align*}
 &\imagpart\left(AY\cdot \Fcred(\omega,Y)\right)\\
 &\qquad=
 -3\omega_0^3 m_1^4 |Y_1|^4
 -3\omega_0^3 m_2^4 |Y_2|^4
 -4\omega_0^3 m_1^2 m_2^2 |Y_1|^2 |Y_2|^2
 +\omega_0^3 m_1^3 m_2 \realpart\left(
 \overline{Y_1}^3 Y_2 - Y_1^3 \overline{Y_2}\right)\\
 &\qquad\le
 -C\omega_0^3 |Y|^4
\end{align*}
for all $(\omega, Y) \in \Hyp \times \C^N$ with 
$$
 A = \left(\begin{matrix}
 m_1^2 &0\\ 0&m_2^2
 \end{matrix}\right).
$$\\
So the condition \eqref{eq_decay_a1} holds
for the system \eqref{eq_dissipation} 
in the case of $m_2 = 3m_1$ with $k=3$.
Therefore by Theorem~\ref{thm_decay}, 
we see that the solution of \eqref{eq_dissipation} 
and its first-order derivatives decay like
$O(t^{-(1/2-1/p)}(\log t)^{-1/2})$ in $L^p$, $p\in [2,\infty]$ as $t\to\infty$.
In this case, however, 
the ratio of masses (i.e. $m_2 = 3m_1$) does not play an essential role
to establish the condition \eqref{eq_decay_a1}.\\

\appendix
\section{Appendix : Proof of Lemma~\ref{lem_energy}}
For the convenience of the readers, 
we give the proof of Lemma~\ref{lem_energy}.
Our goal here is to show $E_s (\tau)\le C\eps^2 \tau^\delta$
under the assumption that $M(T) \le \sqrt \eps$,
where $s \ge 3$ and $0<\delta<1/3$.
Let us define
$$
 E_s (\tau; w_j) 
 = 
 \sum_{k=0}^s \frac{1}{2} \int_{\R} 
 |\pa_\tau \pa_z^k w_j(\tau, z)|^2 
 + \left| \frac{\pa_z}{\tau} \pa_z^k w_j(\tau, z) \right|^2 
 + m_j^2| \pa_z^k w_j(\tau, z)|^2 \,dz
$$
for $s \in \N_0$, $j=1,\ldots,N$, $m_j>0$ 
and a smooth function $w_j$ of $(\tau, z) \in [\tau_0, T) \times \R$.
We start with the following energy inequality:

\begin{lem} \label{lem_energy2}
For $s \in \N_0$, $j=1,\ldots,N$ and $l = 0,1,$ we have
$$
 \frac{d}{d \tau} E_s (\tau;w_j) 
 \le 
 \frac{C'}{\tau^{1+l}}E_{s+l}(\tau;w_j)
 +CE_s(\tau;w_j)^{1/2}\|\Lin_j w_j(\tau,\cdot)\|_{H^s}
 +\frac{C}{\tau^2}E_s(\tau;w_j)
$$
where $L_j$ is given by \eqref{eq_op_lin} and
$C' = 2 \sup_{z \in \R} \frac{|\chi'(z)|}{\chi(z)}$.
\end{lem}

\noindent{\it Proof. }
First, we observe that the following relation holds:
\begin{equation}\label{eq_appendix}
 \left(
  \pa_\tau^2  - \frac{1}{\tau^2}\pa_z^2 +m_j^2\right)\pa_z^k w_j
  =
  \pa_z^k \Lin_j w_j
  +\frac{1}{\tau^{1+l}} \frac{2\chi'(z)}{\chi(z)} 
  \left(\frac{\pa_z}{\tau}\right)^{1-l}\pa_z^{k+l}w_j
  +\frac{1}{\tau^2} P_k(z, \pa_z)w_j
\end{equation}
for any $k\in\N_0$ and $l=0,1$, where
$P_k(z,\partial_z)$ is a $k$-th order operator with bounded coefficients.
Then we can deduce as in the standard energy integral method that
\begin{align*}
 &\frac{d}{d\tau} E_s(\tau;w_j)\\
 &=\sum_{k=0}^s \int_\R (\pa_\tau \pa_z^k w_j)(\pa_\tau^2 \pa_z^k w_j)
   +\frac{1}{\tau^2}(\pa_z^{k+1}w_j)(\pa_\tau \pa_z^{k+1} w_j)
   -\frac{1}{\tau^3}(\pa_z^{k+1}w_j)^2
   +m_j^2(\pa_z^k w_j)(\pa_\tau \pa_z^k w_j)\,dz\\
 &\le \sum_{k=0}^s \int_\R (\pa_\tau\pa_z^k w_j)(\pa_\tau^2 \pa_z^k w_j)
   - \frac{1}{\tau^2} (\pa_z^{k+2}w_j)(\pa_\tau \pa_z^k w_j)
   +m_j^2(\pa_z^k w_j)(\pa_\tau \pa_z^k w_j)\,dz\\
 &\le \sum_{k=0}^s \int_\R |\pa_\tau \pa_z^k w_j| |\pa_z^k \Lin_j w_j|\,dz
   +\frac{C'}{\tau^{1+l}}\sum_{k=0}^s \int_\R |\pa_\tau \pa_z^k w_j|
     \left| \left(\frac{\pa_z}{\tau}\right)^{1-l}\pa_z^{k+l} w_j \right|\,dz\\
 &\qquad+\frac{1}{\tau^2} \sum_{k=0}^s 
    \int_\R |\pa_\tau \pa_z^k w_j||P_k(z,\pa_z)w_j|\,dz\\
 &\le \sum_{k=0}^s \left( \int_\R |\pa_\tau \pa_z^k w_j|^2\,dz\right)^{1/2} 
 \|\pa_z^k \Lin_j w_j\|_{L^2}
   +\frac{C'}{\tau^{1+l}}\sum_{k=0}^s 
   \frac{1}{2}\int_\R |\pa_\tau \pa_z^k w_j|^2
   +\left| \left(\frac{\pa_z}{\tau}\right)^{1-l}\pa_z^{k+l} w_j \right|^2\,dz\\
 &\qquad+\frac{1}{\tau^2} \sum_{k=0}^s \frac{1}{2}
    \int_\R |\pa_\tau \pa_z^k w_j|^2 + |P_k(z,\pa_z)w_j|^2\,dz\\
 &\le CE_s(\tau;w_j)^{1/2}\|\Lin_j w_j\|_{H^s}
   +\frac{C'}{\tau^{1+l}} E_{s+l}(\tau;w_j)
   +\frac{C}{\tau^2}E_s(\tau;w_j),
\end{align*}
where we used \eqref{eq_appendix} for the second inequality. \qed\\ 

We shall apply the above lemma with $l=0$, $s=s_0+s_1+1$ and 
$w_j=v_j$, where $s_0$ is an integer greater than $C'$ 
and $s_1$ is a fixed arbitrary non-negative integer.
Here we regard $G_j$ as a function of 
$(\tau,z,v,\pa_{\tau}v,\pa_z v/\tau)$ for the moment.
Since the Gagliardo-Nirenberg inequality yields 
\begin{align}
 \| 
 \Non_j
   (\tau,\cdot,v(\tau),\pa_\tau v(\tau), \pa_z v(\tau)/\tau)
 \|_{H^s}
 \le 
 \frac{C}{\tau} M(T)^2 E_s(\tau)^{1/2} 
 \le 
 \frac{C}{\tau} \eps E_s(\tau)^{1/2},
 \label{eq_est_G}
\end{align}
we have
$$
 \frac{d}{d\tau} E_{s_0 + s_1 +1}(\tau) 
 \le 
 \left( \frac{C' + C\eps}{\tau} + \frac{C}{\tau^2} \right) E_{s_0 + s_1 +1} (\tau) 
 \le 
 \left(\frac{s_0+\frac{1}{2}}{\tau}
 +\frac{C}{\tau^2}\right)E_{s_0+s_1+1}(\tau)
$$ 
for sufficiently small $\eps$. Thus the Gronwall lemma yields 
$$
 E_{s_0+s_1+1}(\tau) 
 \le 
 E_{s_0+s_1+1}(\tau_0) 
 \exp \left(\int_{\tau_0}^\tau \frac{s_0 + \frac{1}{2}}{\tilde\tau} 
 +\frac{C}{\tilde\tau^2} \, d\tilde\tau \right) 
 \le 
 C \eps^2 \tau^{s_0 + 1/2}.
$$
Next, we apply Lemma~\ref{lem_energy2} with $l=1, s=s_0+s_1$. 
Using the above inequality, we have 
\begin{align*}
 \frac{d}{d \tau} E_{s_0 + s_1} (\tau) 
 &\le
 \frac{C'}{\tau^2} E_{s_0+s_1+1}(\tau) 
 + \frac{C\eps}{\tau}E_{s_0+s_1} (\tau) 
 + \frac{C}{\tau^2}E_{s_0+s_1}(\tau) \\
 &\le 
 C\eps^2 \tau^{s_0-3/2}
 +\left(\frac{C\eps}{\tau}+\frac{C}{\tau^2} \right) E_{s_0+s_1}(\tau).
\end{align*}
Therefore it follows from the Gronwall lemma that 
$E_{s_0+s_1}(\tau) \le C\eps^2 \tau^{s_0 - 1/2}$.
Repeating the same procedure $n$ times, we have 
$E_{s_0+s_1+1 - n}(\tau) \le C\eps^2 \tau^{s_0 - n + 1/2}$
for $n=1,2,\cdots,s_0$. In particular we have 
$$
 E_{s_1 +1} (\tau) 
 \le 
 C \eps^2 \tau^{1/2}.
$$
Finally, we again use Lemma~\ref{lem_energy2} with $l=1, s=s_1$ to obtain 
\begin{align*}
 \frac{d}{d\tau}E_{s_1}(\tau) 
 &\le
 \frac{C'}{\tau^2}E_{s_1+1}(\tau) 
 + \frac{C\eps}{\tau} E_{s_1}(\tau) 
 + \frac{C}{\tau^2}E_{s_1}(\tau)\\
 &\le
 \frac{C\eps^2}{\tau^{3/2}} 
 + \left( \frac{C\eps}{\tau} + \frac{C}{\tau^2} \right) E_{s_1}(\tau).
\end{align*}
The Gronwall lemma yields 
$$
 E_{s_1}(\tau) 
 \le 
 C\eps^2 \tau^{C\eps}
$$
for $\tau \in [\tau_0, T)$. Replacing $s_1$ by $s$ 
and choosing $\eps$ so small that $C\eps \le \delta$, 
we arrive at Lemma~\ref{lem_energy}.\qed\\



\end{document}